\documentclass{amsart}
\usepackage{amssymb,latexsym}
\usepackage{amscd}
\usepackage{hyperref}

\theoremstyle{plain}
\newtheorem{theorem}{Theorem}
\newtheorem{corollary}{Corollary}

\newtheorem{lemma}{Lemma}

\theoremstyle{definition} \newtheorem{definition}{Definition}

\theoremstyle{remark} \newtheorem*{notation}{Notation}

\begin{document}

\title[Construction of Clifford Algebras]{Direct Construction of Grassmann,\\
                   Clifford and Geometric Algebras}
\author{A. Cortzen}

\address{Allan Cortzen}
\email{ac.ga.ca (at) gmail.com}
\urladdr{http://lanco.host22.com/}
\thanks{Thanks goes to the referenced authors.}

\date{2010/11/16}

\begin{abstract}
This is a simple way rigorously to construct Grassmann, Clifford and Geometric Algebras, allowing degenerate bilinear forms, infinite dimension, using fields or certain modules (characteristic 2 with limitation), and characterize the algebras in a coordinate free form.
The construction is done in an orthogonal basis, and the algebras characterized by universality.
The basic properties with short proofs provides a clear foundation for further development of the algebras.
\end{abstract}

\maketitle
\tableofcontents
\section{Introduction}
Often proof of the existence of  Grassmann or Clifford algebras are bypassed, or Chevalleys tensor approach is taken;
 but e.g. investigation of injectivity of the mapping from vector space to Clifford algebra is skipped.
 Pure mathematical books may present a lot of structures before coming to these algebras \hyperlink{r4}{[2,4,9]},
 and the complexity seems to give some problems \hyperlink{r8}{[8]}.\\
Chapters 2, 3 with algebras over $\mathbb{R}^n$ in mind are recommended as minimal reading.
Introductory material can be found in
\hyperlink{r5}{[5]}.

\section{Preliminaries}
Our starting-point may be a field \(R=\mathbb{R}\), a linear space \(V=\mathbb{R}^n\) with basis { }\(e_1=(1,0,\text{$\ldots $0}),e_2=(0,1,\text{$\ldots$0})\ldots \), and indices \(M=\{1,2,\ldots ,\text{\textit{$n$}}\}\) usually ordered by $<$. { }Also used is a quadratic mapping \(q(i)=B(e_i,e_i)\),
where \(B\) is a bilinear form on \(V\) with diagonal form in the basis\( (\left.e_i\right| i\in M)\).
The basic idea behind Grassmann and Clifford algebras is, that products may give new elements, e.g. \(e_1e_2=e_{\{1,2\}}\).\\
In \(\mathbb{R}^n\) the new product should fulfill generator equations \(e_1e_2=-e_2e_1\), \(e_ie_i=q(i)\in \mathbb{R}\) and be associative. Then
a product may be reordered and reduced to get a standard form without repetitions, as in
$$
 e_{\{1,2\}}e_3e_1e_2=-e_1e_2e_1e_3e_2=e_1e_1e_2e_3e_2=-e_1e_1e_2e_2e_3=-q(1)q(2)e_3
 $$
The product properties gives a dimension\(\leq 2^n\), as there are \(2^n\) subsets of \(M\). \\
The goal of the algebra construction is to equip \(W=\mathbb{R}^{2^n}\) with a Clifford product.\\
A basis for \(W\) is \((\left.e_K\right| K\subseteq M)\). Subspaces of \textit{ W} are the scalars \(R e_{\emptyset }\) and \(V\) by identifying
\(e_i\) with \(e_{\{i\}}\).
Sets as indices gives a compact construction. It can indeed be used together with multiindex: \(e_4e_3e_5=e_{(4,3,5)}=-e_{\{4,3,5\}}=-e_{\{3,4,5\}}\).

\begin{notation}
The following notation and definitions will be used. An algebra $A$ is a linear space equipped with a bilinear and associative composition having
a unit \(1_A\). An algebra morphism is supposed to map unit to unit. All the algebras are over the same set of scalars, \(R\). Silently \(x,y\) will
be elements in a linear space \(V\) and \(X,Y\) elements in the algebra at hand.\\
The cardinality of the set \(K \text{ is denoted }|K|\). A product over an index set follow the order given in $M$.
For a index set \(H\) we use \(k<H\) in the meaning \(\forall \, h\in H(k<h)\), implying \(k<\emptyset \) is true.\\
\(H \triangle  J=(H\cup J)\text{$\backslash $(}H\cap J)\) is the symmetric set difference, which is associative.
\end{notation}

\section{Construction }
To make the exposition general we may assume\\
1. \(R\) is a commutative ring with unit \(1\neq 0\)\\
2. \(V\neq \{0\}\) is a free unitary \(R\)-module\\
3. \(B:V\times V\rightarrow R\) is a bilinear form with diagonal form in the basis \( (\left.e_i\,\right | i\in M)\), { }\(q(i)=B(e_i,e_i)\)
and \(M\) is total ordered by a relation \textit{ $<$}.\\
4. { }\(W=\underset{\mathcal{F}}{\oplus } \text{\textit{$R$}}\),
where\textit{  }$\mathcal{F}$ is the set of finite subsets of index set \(M\). \\
A basis for \(W\) is \((\left.e_K\right| K\in \mathcal{F})\). Subspaces of \textit{W} are the scalars \(R\, e_{\emptyset }\) and \(V\) by identifying
\(e_i\) with \(e_{\{i\}}\).\\
\textit{ NB: }If \textit{R} is a field with characteristic different from 2 and \(V\) of finite dimension, any symmetric bilinear form on \(V\)
has an orthogonal basis.\\

Factors $\alpha $ and $\beta $ originate from ordering by swapping and reduction, respectively.
\begin{lemma}
\textit{ Define functions }\(\alpha \)\textit{ and } \(\beta \)\textit{ for combining the sets }\(H,J\in \mathcal{F}\)
\textit{ by}
$$
\alpha (H,J)=\Pi\, (-1) \text{ for } (i,j)\in H\times J \text{ and } j<i
$$
$$
\beta (H,J)=\Pi \, q(i)\text{\textit{$ $}}\textit{ for }\text{\textit{$ $}}i\in H\cap J
$$
\textit{ Then for }\(\psi =\alpha ,\beta \)\textit{  holds}
$$
\psi (H,J)\,\psi (H \triangle  J,K)=\psi (H,J \triangle  K)\,\psi (J,K)
$$
\textit{and therefore also for }\(\psi =\sigma =\alpha \beta \).
\end{lemma}
Proof: As \(\alpha ^2=1\) we get
$$
\alpha (H \triangle  J,K)=\alpha (H\backslash J,K)\,\alpha (J\backslash H,K)\,\alpha (J\cap H,K)^2=\alpha
(H,K)\,\alpha (J,K)
$$
and likewise \(\alpha (H,J \triangle  K)=\alpha (H,J)\,\alpha (H,K)\). Thus the equation for $\alpha$  is obvious.\\
By Venn diagrams
$$\beta (H,J)\,\beta (H \triangle  J,K)= \Pi\, q(i)
 \text{ for \(i\) in just two of the sets \(H, J, K\), }
$$
 and \(\beta (H,J \triangle  K)\beta (J,K)\) gives the same result.

\begin{theorem}
Define a product \((X,Y)\rightarrow X Y\)\textit{  in W by} \(e_H \, e_J=\sigma (H,J)\,e_{H \triangle  J}\)\textit{  and bilinearity. Then W becomes an algebra with }
\(e_{\emptyset }\)\textit{  as unit, such that}
$$
i\neq j\Rightarrow e_ie_j=-e_je_i\textit{ and }x\in
V\Rightarrow x^2=B(x,x)\,e_{\emptyset }
$$
\end{theorem}
Proof: The associative law \((x y) z=x (y z)\) is multilinear and is verified for basis elements by
$$
(e_H e_J)\,e_K=\sigma (H,J)\,e_{H \triangle  J}\, e_K=\sigma (H,J)\,\sigma (H \triangle  J,K)\,e_{(H
\triangle  J)\triangle  K}
$$
$$
e_H \,(e_Je_K)=\sigma (J,K)\,e_He_{J \triangle  K}=\sigma (H,J \triangle  K)\,\sigma (J,K)\,e_{H
\triangle  (J \triangle  K)}
$$
Furthermore
$$
e_{\emptyset }e_K=\sigma (\emptyset ,K)e_K=e_K=\sigma (K,\emptyset )e_K=e_Ke_{\emptyset }\textit{ and } e_ie_i=\sigma
(\{i\},\{i\})e_{\emptyset } =q(i)\,e_{\emptyset }.
$$
Now \(i<j\Rightarrow e_ie_j=\sigma (\{i\},\{j\})\,e_{\{i,j\}} =e_{\{i,j\}}\) and likewise \(e_je_i=-e_{\{i,j\}}\).
If \(x=\Sigma _i\,\lambda _ie_i\), then by separating into cases $i<j \,, i=j$ { and }$j<i$, we find
$$
x^2=\Sigma _{i,j}\lambda _i\lambda _je_ie_j=\Sigma  \lambda _i^2e_i^2e_{\emptyset
}=B(x,x)e_{\emptyset }.
$$
\begin{corollary}
\(e_K=\Pi  _{i\in K} e_i\)
\end{corollary}
Proof: Use induction and \(j<H\Rightarrow e_je_H=\sigma (\{j\},H)e_{\{j\}\cup H} =e_{\{j\}\cup H}\).

\begin{corollary}
For any algebra \(A\) over F and any linear mapping \(f:V\rightarrow A\) with
$$f(x)^2=B(x,x)\,1_A$$
there exists a unique algebra morphism \(F:\text{\textit{$ $}}W\rightarrow A\), which extends f.
\end{corollary}

Proof: Define a linear mapping \(F:\text{\textit{$ $}}W\rightarrow A\) necessarily by
$$
F(e_{\emptyset })=1_A \text{ and } F(e_K)=\Pi  _{i\in K} f(e_i)
$$
For \(i\neq j\) and \(x=e_i+e_j\), we get
$$
f(x)^2=B(x,x)1_A\Rightarrow f(e_i)f(e_j)=-f(e_j)f(e_i).
$$
Let \(H=\left\{h_1,h_2,\ldots ,h_p\right\}\) and \(K=\left\{k_1,\ldots ,k_q\right\}\) with \(h_i<h_{i+1}\) and  \(k_j<k_{j+1}\)
for relevant indices. Ordering from lower to higher indices by swapping neighbors in
$$
F(e_H) F(e_J)=f(e_{h_1}) f(e_{h_2})\ldots  f(e_{h_p})f(e_{k_1})f(e_{k_2})\ldots
 f(e_{k_q})
$$
 gives a sign change \(\alpha (H,K)\), and reducing equals gives a further factor \(\beta (H,K)\). Thus \(F(e_H) F(e_J)=\sigma (H,J)F(e_{H \triangle  J})=F(e_He_J)\).\\

\begin{definition}
\textit{A Clifford algebra over $B$, denoted $\textit{Cl}_V(B)$ or $\textit{Cl}(B),$
is an algebra isomorphic to \(W_B=W\) with an isomorphism fixing \(V\).
In the case \(B=0\), we have a Grassmann algebra, \(\Lambda (V)=\textit{Cl}_V(0)\), and a product, the outer product }$\wedge$.
\end{definition}

\begin{definition}
\textit{If the linear space of the Clifford algebra }\(W\)\textit{ also is given a Grassmann structure by the zero bilinear form, we get a double algebra }\(W_{B,0}=W\).\\
\textit{A geometric algebra over }\(B\)\textit{ , denoted }\(\Lambda (V,B)\),\textit{  is a double algebra isomorphic
to }\(W_{B,0}\)\textit{  with an isomorphism fixing }\(V\).
\end{definition}

\begin{corollary}
\textit{In }\(\Lambda (V)\)\textit{  define the subspace of elements of grade }\(r\in \mathbb{Z}\)\textit{  by\\
}\textit{  }\(\Lambda _r(V)=\text{span}\left\{\wedge _{i=1}^ra_i|a_i\in V\right\}\)\textit{  for }\(r\geq 0\)\textit{  and otherwise}\textit{  }\(\Lambda
_r(V)=\{0\}\)\textit{ . Then\\
1. }\textit{  }\(\Lambda (V)=\oplus _r\text{\textit{$ $}}\Lambda _r(V)\)\textit{  and}\textit{  }\(\Lambda _r(V)\wedge \Lambda _s(V)\subseteq \text{\textit{$
$}}\Lambda _{r+s}(V)\)\textit{ .\\
This allows us to define }\(X\rightarrow \langle X\rangle _r\)\textit{ to be the projection on}\textit{  }\(\text{\textit{$ $}}\Lambda _r(V)\)\textit{  along }\(\oplus _{i\neq
r}\text{\textit{$ $}}\Lambda _i(V)\)\\
\textit{ 2. { }}\(x_1\wedge x_2=-x_2\wedge x_1\)\textit{ \\
3. { }}\(x_1\wedge x_2\wedge \ldots  \wedge x_p\)\textit{  { }is multilinear and alternating in the x-variables\\
4. { }}\(\left(x_1,x_2\ldots  ,x_p\right)\)\textit{  is linear independent { }}$\Leftrightarrow $\textit{  }\(x_1\wedge x_2\wedge \ldots  \wedge
x_p\) \textit{  is linear independent. }\\
5. For any algebra \(A\)\textit{ over R and linear mapping }\(f:V\rightarrow A\)\textit{ with }\(f(x)^2=0\), \textit{there exists a unique algebra morphism or outermorphism }\(f_{\wedge }:\Lambda (V)\rightarrow A\)
\textit{extending f.}
\end{corollary}

Proof: 1. Observe that 0 can be assigned any appropriate grade.
In\textit{  }\(W_0\) define \(\Omega _r=\text{span}\left\{\left.e_K\right||K|=r,K\in \mathcal{F}\right\}\) for \(r\geq 0\) and otherwise \(\Omega _r=\{0\}\)\textit{
.} Then clearly
\(W_0=\oplus _r\text{\textit{$ $}}\Omega _r\) and \(\Omega _r\wedge \Omega _s\subseteq \text{\textit{$ $}}\Omega _{r+s}\), which implies\textit{
 }\(\wedge _{i=1}^ra_i\in \text{\textit{$ $}}\Omega _r\). Thus\textit{  }\(\Omega _r=\Lambda _r\). \\
Now use the isomorphism \(W_0\rightarrow \Lambda (V)\) fixing \(V\).\\
2. In\textit{  { }}\(x\wedge x=0\) set \textit{  }\(x=x_1+x_2\). \\
3. Obviously the expression is multilinear, and { }0, if\textit{  }\(x_i=x_{i+1}\). For  \(i<j\) and \textit{  }\(x_i=x_j\) this situation can
be obtained by swapping neighbors.\\
4. The full proof is in the appendix.\\
5. Consequence of corollary 2.

\begin{corollary}
\textit{
In }\(\text{\textit{$\text{\textit{Cl}}_V$}}(B)\)\textit{  define the subspace of grade }\(r\in \mathbb{Z}_2=\{0,1\}\) \textit{ by}
$$\text{\textit{Cl}}_{V,r}(B)=\text{span}\left\{\Pi _{i=1}^sa_i \,|s\equiv r\text{  }(\text{mod } 2), a_i\in V\right\}$$

\textit{Then}\textit{  }\(\text{\textit{Cl}}_V(B)=\oplus _r \text{\textit{Cl}}_{V,r}(B)\)\textit{  and}\textit{  }\(\text{\textit{Cl}}_{V,r}(B) \text{\textit{Cl}}_{V,s}(B)\subseteq
\text{\textit{Cl}}_{V,r+s}(B)\).
\end{corollary}

Proof: In\textit{  }\(W_B\) define \(\Omega _r=\text{span}\left\{\left.e_K\right||K| \equiv r\text{  }(\text{mod } 2),K\in \mathcal{F}\right\}\)\textit{ .} Then we have\\
\(W_B=\oplus _r\text{\textit{$ $}}\Omega _r\) and \(\Omega _r\Omega _s\subseteq \text{\textit{$ $}}\Omega _{r+s}\), since factor reductions for products
are even in number. This implies\textit{  }\(\Pi _{i=1}^sa_i\in \text{\textit{$ $}}\Omega _r\) for \(s \equiv r\text{  }(\text{mod } 2) \). Hence\textit{  }\(\Omega _r=\text{\textit{Cl}}_{V,r}\).\\
Now use the isomorphism \(W_B\rightarrow \text{\textit{$\text{\textit{Cl}}_V$}}(B)\) fixing \(V\).

\section{Characterization of universal Clifford algebras}

\begin{theorem}
\textit{Let }\(\mathcal{A}_R(V,B)\)\textit{  be the category of linear mappings }\(f\)\textit{  from V into an algebra A, such that} \(f(x)^2= B(x,x)1_A\)\textit{
. \\
A mapping }\(\rho :V\rightarrow U\)\textit{  in }\(\mathcal{A}_R(V,B)\)\textit{  is said to be universal, if for every linear mapping }\(f:V\rightarrow
A\)\textit{  in }\(\mathcal{A}_R(V,B)\), \textit{  there is a unique algebra morphism }\(F:\text{\textit{$ $}}U\rightarrow
A\)\textit{  such that }\(F\circ \rho =f\).
\textit{ \\
In case} \(V\subset U\)\textit{  this means }\(F\)\textit{  extends f , when }\(\rho \)\textit{  silently is taken as the injection.} \\
\textit{If { }}\(f_i:V\rightarrow U_i, i=1,2\)\textit{  are universal in }\(\mathcal{A}_R(V,B)\),\textit{  then there exists a unique algebra isomorphism
}\(F:\text{\textit{$ $}}U_1\rightarrow U_2\)\textit{ such that }\(F\circ f_1=f_2\).
\end{theorem}
Proof:\textit{  }Universality gives unique algebra morphisms\textit{  }\(F:\text{\textit{$ $}}U_1\rightarrow U_2\), \textit{  }\(G:\text{\textit{$
$}}U_2\rightarrow U_1\), such that \(F\circ f_1=f_2\), { }\(G\circ f_2=f_1\).
As { }\(F\circ G\circ f_2=f_2\) and { }\(\text{id}_{U_2}\circ f_2=f_2\), universality implies \(\text{id}_{U_2}=F\circ G\), and likewise \(\text{id}_{U_1}=G\circ
F\).
\begin{corollary}
\(\text{\textit{$\text{\textit{Cl}}_V$}}(B)\)\textit{  is universal in }\(\mathcal{A}_R(V,B)\).
\end{corollary}

Proof: Follows from corollary 2 and definition 1.

\begin{corollary}
The Clifford product in
\(\text{\textit{$\text{\textit{Cl}}_V$}}(B)\) is independent of the orthogonal basis in theorem 1.
\end{corollary}

Proof: Let \textit{  }\(\left(\ddot{e}_i|i\in \ddot{M}\right)\) be an orthogonal basis for \(V\), and \(\ddot{W}_B\) the Clifford algebra constructed as in theorem 1.
By universality we now get an unique isomorphism \(\ddot{F}:\ddot{W}_B\rightarrow \text{\textit{$\text{\textit{Cl}}_V$}}(B)\) fixing \(V\).
Hence \(\text{\textit{Cl}}_V(B)\) has the product defined from \(\ddot{W}_B\).

\begin{corollary}
In \(\text{\textit{$\text{\textit{Cl}}_V$}}(B)\)\textit{  the main automorphism }\(X\rightarrow \hat{X}\)\textit{  is the universal extension
of }\(f:V\rightarrow \text{\textit{$\text{\textit{Cl}}_V$}}(B)\)\textit{  where { }}\(f(x)=-x\).\\
 \textit{ We have }\(X\rightarrow \hat{X}=(-1)^rX\)\textit{
 for X of Clifford grade r}.
\end{corollary}

Proof: By universality\(f\) can be extended { }uniquely to an algebra morphism \(\hat{x}:\text{\textit{$ $}}\text{\textit{$\text{\textit{Cl}}_V$}}(B)\rightarrow
\text{\textit{$\text{\textit{Cl}}_V$}}(B)\). In a basis by linearity and \(\hat{e_K}=\Pi  _{i\in K} f\left(e_i\right)=(-1)^{|K|}e_K=(-1)^{|K| \bmod
2}e_K\) { }the statement is proved.

\begin{corollary}
In $Cl_V(B)$ \textit{define the reversion }\(\widetilde{x}\)
\textit{ by }
\((X Y)^{\sim }=\widetilde{Y}\widetilde{X}\)\textit{, linearity and by fixing }$1_{Cl(B)}$\textit{ and }\(V\). Then \(X^{\sim \sim }=X\)\textit{, and { }}\(\left(a_1a_2\ldots  a_r\right){}^{\sim }=a_r\ldots  a_2 a_1\).
\end{corollary}

Proof: In the linear space $U$ of $Cl_V(B)$ an algebra \((U,\diamond )\) is defined by the product \(X\diamond
Y=Y X\). As { }\textit{  }\(x\diamond x= B(x,x)1_U\) and
$$
(X Y)^{\sim }=\widetilde{Y} \widetilde{X}\Leftrightarrow (X Y)^{\sim }=\widetilde{X}\diamond \widetilde{Y}
$$
a reversion must be an algebra morphism \(\text{\textit{  }}\widetilde{x}:Cl_V(B)\rightarrow (U,\diamond )\), and universality of $Cl_V(B)$
implies uniqueness and existence. Now \((X Y)^{\sim \sim }=\left(\widetilde{Y} \widetilde{X}\right)^{\sim }=X Y\) and by universality\textit{  }\(X^{\sim
\sim }=X\).
The last formula follows from \((X Y)^{\sim }=\widetilde{Y}
\widetilde{X}\).

\section{Characterization of Geometric algebras}
\begin{definition}
Set \(\chi (S)=1\)\textit{, if }\(S\)\textit{ is true, and else zero}.\\
\textit{ In the geometric algebra }\(W_{B,0}\) \textit{
define mappings} \(*,\rfloor \text{\textit{$ $}}\text{\textit{and}}\text{\textit{$ $}}\lfloor \)\textit{ by bilinearity and \\
}\(e_H *e_K=\chi (H=K) \, e_H e_K\) (\textit{the scalar product}),\\
\(\left.e_H \right\rfloor e_K=\chi (H\subseteq K) \, e_H e_K\)   (\textit{the left contraction}),\\
\(e_H \left\lfloor e_K=\chi (H\supseteq K) \, e_H e_K\right.\)   (\textit{the right contraction}).\\
\textit{ Observe that in a geometric algebra the grading is taken from the Grassmann structure, and that }\(e_H \wedge e_K=\chi (H\cap K=\emptyset
) \, e_H e_K\).
\end{definition}

\begin{theorem}
In a geometric algebra \(\Lambda (V,B)\)\textit{  holds\\
1. }\(1_{\text{\textit{$\text{Cl}(B)$}}}=1_{\Lambda }\)\textit{ and} \(x X=x\wedge X+x\rfloor X\)\\
2. The two subalgebras have the same reversions, and $$\widetilde{X}=(-1)^{r(r-1)/2}X \textit{  for X of grade r}$$
3a. $
(X \lfloor Y )^{\sim }=\widetilde{Y} \rfloor \widetilde{X},\textit{  }
X*Y=\langle X Y\rangle _0,\textit{ and }
X\wedge Y=\langle X Y\rangle _{r+s},
$\\
3b. If \(\text{grade}(X)=r\)\textit{  and { }}\(\text{grade}(Y)=s\)\textit{, then}
$$
X\rfloor Y=\langle X Y\rangle _{s-r}
\textit{, }
X\lfloor Y=\langle X Y\rangle _{r-s}
\textit{ and }
X Y=\Sigma
_{i=|r-s|\text{, step}\text{  }2}^{r+s}\langle X Y\rangle _i
$$
4. \((X\wedge Y )\rfloor Z=X\rfloor (Y \rfloor Z)\)\\
5. \(x \rfloor y=B(x,y)1_{\Lambda }\) \textit{ and }\(x \rfloor (X Y)=(x \rfloor X) Y+\hat{X} (x \rfloor Y)\)\\
6. \(x \rfloor y=B(x,y)1_{\Lambda }\) \textit{ and }\(x \rfloor (X\wedge Y)=(x \rfloor X)\wedge Y+\hat{X}\wedge (x \rfloor Y)\)\\
7. \(x\rfloor (x_1x_2\ldots  x_p)=\sum _{ k=1}^p (-1)^{k-1}x_1x_2\ldots  (x\rfloor x_k)\ldots  x_p\)\\
8. \(x\rfloor (x_1\wedge x_2\wedge \ldots  \wedge x_p)=\sum _{ k=1}^p (-1)^{k-1}x_1\wedge x_2\ldots \wedge  (x\rfloor x_k)\ldots \wedge x_p\)\\
9. \( x_1,x_2,\ldots  ,x_p\)\textit{ are pairwise orthogonal }$\Rightarrow $
\textit{ }\(\Pi _{i=1}^p x_i=\wedge _{i=1}^p x_i\)
\end{theorem}

Proof: By isomorphy and linearity it should be sufficient to sketch proofs  of the statements for basis elements in \(W_{B,0}\).\\
1. Now \(\left.h\in H\Rightarrow e_h \wedge e_H+e_h \right\rfloor
e_H=0+e_h\text{  }e_H\) and similar for \(h\notin H\). \\
2. To reorder \(\widetilde{e_H}\) requires \((-1)^{|H|(|H|-1)/2}\) swappings in boths structures.\\
3. The proofs are much alike, so we will take some examples:\\
The first in (3a): \(\left(e_H \left\lfloor e_J\right){}^{\sim }=\chi (J\subseteq H)\left(e_H e_J\right){}^{\sim }=\chi (J\subseteq
H)\widetilde{e_J} \widetilde{e_H}=\widetilde{e_J} \right\rfloor \widetilde{e_H}\).\\
The last in (3b) follows from \(e_H e_J=\sigma (H,J)e_{H \triangle  J}\) and\\ \(|H \triangle  J|=|H|+|J|-2|H\cap J|=|H|-|J|+2|J\backslash H|=|J|-|H|+2|H\backslash J|\)\\
4. Can be reduced to \(\chi (H\cap J=\emptyset )\,\chi (H,J\subseteq K)=\chi (H\subseteq (K\backslash J))\,\chi (J\subseteq K)\)\\
5. Clearly \(\left.e_i \right\rfloor e_j=B\left(e_i,e_j\right)e_{\emptyset }\) and the remaining can be reduced to
$$
\chi (i\in (H \triangle  J))=\chi (i\in H)+(-\chi (i\in H)\,\chi (i\in J)+\chi (i\notin H)\,\chi (i\in J))
$$
6. For \(\text{grade}(X)=r\) and \(\text{grade}(Y)=s\) use (3) and take grade \(r+s-1\) in (5)\\
7,\,8. Follows from (5,\,6), as e.g.
$$\left.\left.x\rfloor \left(x_1\left(x_2\ldots  x_p\right)\right)=(x \rfloor x_1\right)\text{ }\left(x_2\ldots
 x_p\right)-x_1(x \rfloor  \left(x_2\ldots  x_p\right)\right),
$$
$$
x \rfloor(x_2\ldots  x_p)=(x \rfloor x_2)(x_2\ldots  x_p)-x_2(x \rfloor  (x_2\ldots  x_p)) \textnormal{ etc.}
$$
9. By (1,7), as e.g.
$$ x_1 x_2\ldots  x_p=x_1\wedge (x_2\ldots  x_p)+x_1\rfloor  (x_2\ldots  x_p)=x_1\wedge
 (x_2\ldots  x_p) \textnormal{ etc.}
$$

\begin{theorem}
Given a double algebra U of \(\Lambda (V)\)\textit{ and }\(\text{\textit{$\text{\textit{Cl}}_V$}}(B)\)\textit{  occupying the same linear space, such that }\(1_{\text{\textit{Cl}}(B)}=1_{\Lambda }\)\textit{, and V is a common linear space. Then}
$$
 U \textit{  is a geometric algebra over B}
$$
$$
 \Leftrightarrow \textit{ } \forall a_1,\ldots  ,a_p\in V\left( a_1,a_2,\ldots  ,a_p\right.
\textit{ are pairwise orthogonal}
\Rightarrow  \Pi _{i=1}^p a_i=\wedge _{i=1}^p a_i).
$$
\end{theorem}

Proof: $\Rightarrow $: Follows from lemma 1.\\
$\Leftarrow $: By universality let \(\phi :\Lambda (V,B)\rightarrow U \) be the unique isomorphism fixing \(V\) determined be the Grassmann structures and similar \(\psi :\Lambda (V,B)\rightarrow U \) for the Clifford structures.
Let \(\left(\left.e_i\right| i\in M\right)\) be an orthogonal basis for \(V\).
In \(\Lambda (V,B\)) set \(\text{\textit{$ $}}X=\wedge _{i\in K} e_i=\Pi _{i\in K} e_i\).
Then \(\phi =\psi \), as
$$
\phi ( X)=\wedge _{i\in K}\phi \left( e_i\right)=\Pi _{i\in K}\psi \left(e_i\right)=\psi (X).
$$

\begin{corollary}
If \(U\) is a geometric algebra over V, then the Clifford structure determines the Grassmann structure and B
 - and vice versa.
\end{corollary}

Proof: By universality let \(\phi :\Lambda (V,B)\rightarrow U\) be the unique isomorphism fixing \(V\) determined by the Clifford structures.
Then $\phi$ determines the Grassmann structure of \(U\). The converse is similar proven.\\

In \hyperlink{r3}{[3]} various Grassmann structures in a given Clifford algebra are used to describe interacting fermions.

\section{Construction of Clifford algebras from tensor algebras}
For the tensor algebra over \(V\), \(\mathcal{T}=\mathcal{T}(V,\otimes )\)\text{, this universality statement is valid}:\\
 \textit{ For any algebra
A over R and any linear mapping }\(g:V\rightarrow A\)\textit{ , there is an unique algebra morphism }\(G:\mathcal{T}\rightarrow A\)\textit{  which
extends g\\
}\\
Suppose $D$ is any bilinear form on $V$. We may then extend theorem 2 by replacing \(B\) with \(D\) and in this way extend the Clifford algebra concept.

\begin{definition}
\textit{A Clifford algebra over \(D\) is an universal object in \(\mathcal{A}_R(V,D)\)}
\end{definition}

\begin{theorem}
Let \(\mathcal{I}=\mathcal{I}(V,D)\)\textit{ be the two-sided ideal in} \(\mathcal{T}=\mathcal{T}(V)\)
\textit{ generated by }\(\mathcal{S}=\{x\otimes
x-D(x,x)1_{\mathcal{T}}|x\in V\}\)\textit{, }\(\text{\textit{Cl}}=\mathcal{T}/\mathcal{I}\)\textit{ the quotient algebra and }\(\hat{\pi}:\mathcal{T}\rightarrow \text{\textit{Cl}}\)
\textit{ the quotient mapping. Then }\(\pi =\hat{\pi}|_V\)
\textit{is an universal object in }\(\mathcal{A}_R(V,D)\).\\
\textit{If }\(D=B\)\textit{, then }\(\pi\)\textit{ is injective.}
\end{theorem}
\[
\begin{CD}
V @>\subset>> \mathcal{T} @>\hat{\pi}>>\mathcal{T}/\mathcal{I}=Cl\\
@V id VV @V G VV @VV F V \\
V @> f >> A @> id >> A
\end{CD}
\]
\\

Proof: In \(\text{\textit{Cl}}\) elements are of the form \(t+\mathcal{I},\, t\in \mathcal{T}\),
and obviously \(x\otimes x+\mathcal{I}=D(x,x)1_{\mathcal{T}}+\mathcal{I}\).
Therefore \(\pi (x)^2=D(x,x)1_{\text{\textit{Cl}}}\) and \(\pi\) is an object in \(\mathcal{A}_R(V,D)\).\\
To prove universality of \(\pi\) we will use universality for tensors.
Therefore, to any object \(f:V\rightarrow A\) in \(\mathcal{A}_R(V,D)\) there is a unique algebra morphism \(G:\mathcal{T}\rightarrow A\) which extends
\textit{f}. As \(f(x)^2=D(x,x)1_A\) implies \(\mathcal{S}\subseteq G^{-1}(0)\), there exists a unique algebra morphism \(F:\text{\textit{Cl}}\rightarrow
A\), such that \(G=F\circ \hat{\pi }\). Hence \(f=F\circ \pi \). \\
Conversely, as any algebra morphism \(F:C\text{\textit{$l$}}\rightarrow A\), gives an algebra morphism \(F\circ \hat{\pi }\) extending
\(f\), we have \(G=F\circ \hat{\pi }\) and \(F\) is unique. Thus \(\text{\textit{Cl}}\) is universal in \(\mathcal{A}_R(V,D)\). \\
If \(D=B\), by universality of \(\text{\textit{Cl}}\) there exists an algebra morphism \(F:\text{\textit{$ $}}C\text{\textit{$l$}}\rightarrow
\text{\textit{$\text{\textit{Cl}}_V$}}(B)\), 
such that \(F\circ \pi =\text{id}_V\), and therefore \(\pi\) is injective.

\section{Conclusion}
On elementary basis we have defined and constructed the different algebras.\\
The universality principle has been described, and used in many ways:
\begin{itemize}
\item To prove in-dependency of orthogonal basis.
\item To define the main automorphism and the reversion.
\item In various proofs.
\item To fully define Clifford algebras
\item To establish connection to Chevalley's tensor based construction.
\end{itemize}
Comments on other constructions of geometric algebra can be found in \hyperlink{r7}{[7]}.\\
A more general construction of Clifford algebras over modules is found in \hyperlink{r4}{[2,4]}, and in a forthcoming paper.

\section{Appendix}

\begin{theorem}
In \(\Lambda (V)\)\ holds
$$
S=(x_1,x_2\ldots  ,x_p) \textit{ is linear dependent }\Leftrightarrow
S_{\wedge}=\wedge _{k=1}^px_k \textit{ is linear dependent.}
$$
\end{theorem}

Proof: $\Rightarrow $: Assume \(\Sigma  \lambda _ix_i=0\) and \(\lambda _j\neq 0\), then \(\lambda _jx_j=-\Sigma _{i\neq j} \lambda _ix_i=0\), which
by corollary 3 (3) gives \(\lambda _jS_{\wedge }=0\).\\
$\Leftarrow $: If \(R\) is a field, assume \(S\) is linear independent and
 construct \(\Lambda (V)\) from a basis containing \(S\) (corollary
6). As \(S_{\wedge}\) is a basis element, \(\lambda \neq 0\Rightarrow \lambda  S_{\wedge}\neq 0\).\\
$\Leftarrow $: Obvious for \(p=1\). Assume that \(S\) is linear independent, and \(S_{\wedge}\) is linear dependent and \(p\) is the smallest
number, for which such a set \(S\) can be found. Hence, if \(T\) is a strict subset of \(S\),
then \(T_{\wedge}=\wedge _{x \in T}x\) is linear independent.\\
Let \(\left(\left.e_i\right|i\in M\right)\) be a basis for \(V\) and define a geometric algebra structure on
\(\Lambda (V)\) by
letting \(\left(e_i\right)\) be an orthogonal basis and \(q\left(e_i\right)=1\).
Set \(X= x_1\wedge x_2\wedge \ldots  \wedge x_{p-1}\). Then \(\lambda  X\wedge x_p=0\) for some \(\lambda \neq 0\).
Furthermore \(\lambda  X\neq 0\), \(p\geq 2\) and \(\lambda  x_p\neq 0\).\\
From \(0\neq \lambda  X=\Sigma _{K\in \mathcal{E}} \lambda _Ke_K\) select \(\lambda _K\neq 0\).
By theorem 3 (8)
$$
e_{i_{p-1}}\rfloor (\ldots (e_{i_1}\rfloor (x_1\wedge \ldots  \wedge x_p)))=\Sigma _{ k=1}^p \mu _k x_k,
$$
as \(p-1\) of { }the \(x_i\)-elements are contracted. Hence\\
$$
0=\widetilde{e_K}\rfloor (\lambda  X\wedge x_p) =\lambda _Kx_p+\Sigma _{k=1}^{p-1} \mu _kx_k,
$$
which contradicts the assumed in-dependency of \(S\).

\nonumber

\begin{gather}
\end{gather}

\end{document}